\documentclass[preprint,review,12pt,3p]{elsarticle}
\usepackage{amsmath}
\usepackage{graphicx,amssymb,psfrag,subfigure,color}
\usepackage{enumitem}
\usepackage{epsfig}
\usepackage{bm}
\usepackage{cases}
\usepackage{upgreek}
\usepackage{array}
\usepackage{lineno}
\usepackage[colorlinks]{hyperref}
\usepackage{epstopdf}
\usepackage{booktabs}
\usepackage{natbib}
\usepackage[table]{xcolor}
\usepackage[section]{placeins}
\usepackage{multirow}
\usepackage{makecell}
\usepackage{diagbox}

\makeatletter
\renewcommand{\maketag@@@}[1]{\hbox{\m@th\normalsize\normalfont#1}}%
\makeatother
\newcommand{\beq}{\begin{equation}}
\newcommand{\eeq}{\end{equation}}
\newcommand{\bpm}{\begin{pmatrix}}
\newcommand{\epm}{\end{pmatrix}}
\newcommand{\beqa}{\begin{eqnarray}}
\newcommand{\eeqa}{\end{eqnarray}}
\newcommand{\beqas}{\begin{eqnarray*}}
\newcommand{\eeqas}{\end{eqnarray*}}

\newcommand{\non}{\nonumber}

\def\XXint#1#2#3{{\setbox0=\hbox{$#1{#2#3}{\int}$ }
\vcenter{\hbox{$#2#3$ }}\kern-.6\wd0}}

\biboptions{sort&compress}
\begin{document}

\begin{frontmatter}

\title{Stiffness minimisation of graded  microstructural configurations using asymptotic analysis and machine learning}

%\tnotetext[mytitlenote]{Fully documented templates are available in the elsarticle package on \href{http://www.ctan.org/tex-archive/macros/latex/contrib/elsarticle}{CTAN}.}

%% Group authors per affiliation:
\author[DUT1,DUT2]{Chuang Ma}
\author[PSU]{Dingchuan Xue}
\author[DUT1,DUT2]{Shaoshuai Li}
\author[DUT1,DUT2]{Zhengcheng Zhou}
\author[DUT1,DUT2,DUT3]{Yichao Zhu\corref{cor1}}
\ead{yichaozhu@dlut.edu.cn}
\author[DUT1,DUT2,DUT3]{Xu Guo\corref{cor1}}
\ead{guoxu@dlut.edu.cn}
\cortext[cor1]{Corresponding author}
\address[DUT1]{State Key Laboratory of Structural Analysis for Industrial Equipment, Department of Engineering Mechanics, Dalian University of Technology, Dalian, 116023, P. R. China}
\address[DUT2]{International Research Center for Computational Mechanics, Dalian University of Technology}
\address[PSU]{Department of Engineering Science and Mechanics, Pennsylvania State University, University Park, PA 16802, USA}
\address[DUT3]{Ningbo Institute of Dalian University of Technology, No.26 Yucai Road, Jiangbei District, Ningbo, 315016, P. R. China}
%\fntext[myfootnote]{Since 1880.}

\begin{abstract}
  The article is aimed to address a mutually boosting use of asymptotic analysis and machine learning, for fast stiffness design of configurations infilled with smoothly-varying graded microstructures. The discussion is conducted in the context of an improved asymptotic-homogenisation topology optimisation (AHTO plus) framework \cite{zhu2019novel}. It is demonstrated that on one hand, machine learning can be employed to represent the key but implicit inter-relationships revealed from asymptotic analysis, and the evaluations of the homogenised quantities, as well as the sensitivities of the design variables, become quite efficient. On the other hand, the use of asymptotic analysis identifies a computational routine for data acquisition, thus the training data here are inexhaustible in theory. Key issues regarding integration of the two methods, such as ensuring the positive definiteness of the homogenised elasticity tensor represented with neural networks, are also discussed. The accuracies and the efficiencies of the present scheme are numerically demonstrated. For two-dimensional optimisation, it takes the present algorithm roughly 300 seconds on a standard desktop computer, and this qualifies the present scheme as one of the most efficient algorithms used for the compliance optimisation of configurations infilled with complex microstructures.
\end{abstract}

\begin{keyword}
  Graded microstructural configuration \sep Multiscale \sep Machine learning\sep Asymptotic analysis \sep Optimal design of structural compliance
% \MSC[2010] 00-01\sep  99-00
\end{keyword}
\end{frontmatter}

%\linenumbers

\section{Introduction}
Configurations infilled with lattice cells have been demonstrated to be perspective for high-end applications related to mechanics, acoustics and optics, etc. \cite{lakes1993materials, sigmund1996composites, kushwaha1993acoustic, liu2015frequency, aage2017giga}. Recently with the rapid development of additive manufacturing technologies, the fabrication techniques of functional components and devices decorated with finely designed microstructures have become more and more mature, and this naturally promotes the development of corresponding intelligent/automatic design algorithms. This normally entails the use of multiscale approaches, where a key issue is to manage the balance between the quality of results and computational efficiency \cite{WuJ_SMO_Review2021}.

When the constituting unit of a lattice configuration is periodic in space, its certain mechanical behaviours can be effectively analysed with the asymptotic homogenisation method initiated by Bensoussan et al. \cite{papanicolau1978asymptotic}. The key idea is to treat the lattice configuration as a homogeneous continuum, and the equivalent elastic moduli are then computed by solving a set of cell problems defined within one (spatially periodic) cell. The method was then employed for the optimal design of spatially periodic configurations by Bendsoe and Kichuchi \cite{bendsoe1989optimal}, where the constituting cells are permitted to (uniformly) change their size and rotational angles. Studies in the field were then advanced with the proposal of the so-called inverse homogenisation method \cite{sigmund1994materials}, where the structure of the constituting cells is tuned so as to meet certain requirements over the system performance. A number of succeeding studies were conducted, seeking for appropriate microstructural design of unit cells to fulfil multi-functional purpose, such as the reliable service in extreme performance materials \cite{huang2012evolutionary, zhou2012design, radman2013topological}, piezoelectric material microstructure design \cite{silva1997optimal}, band gap material design \cite{sigmund2003systematic}, material permeability optimisation \cite{guest2007design}, etc. Studies were also extended to consider optimal design on concurrent length scales, where the design variables originate from both the microscopic topology and the macroscopic material density distribution \cite{rodrigues2002hierarchical, liu2008optimum, coelho2008hierarchical}. The concurrent optimal design was also considered for other engineering purpose, such as multi-functional design \cite{niu2009optimum,yan2016multi,deng2013multi} and uncertain loading design \cite{deng2017concurrent}, etc.

It is noted that the asymptotic homogenisation formulation underpinning the aforementioned works was actually derived for configurations where the infilled microstructures are periodic in space. However, a large number of perspective microstructural configurations, either naturally formed \cite{sanchez2005biomimetism, fratzl2009biomaterial, meyers2013structural} or artificially made \cite{jorgensen1998spherical, arabnejad2012multiscale, cheng2018extra}, bear heterogeneous microstructures. Then the presumption about microstructural periodicity breaks down, and the asymptotic homogenisation theory needs modification, so as to fuel the development of the design methods for graded microstructural configurations (GMCs).

For this purpose, several research works \cite{zhou2008design, wang2017concurrent, zhang2018multiscale, zhang2019concurrent} were proposed. Given the improvements brought by them, the above-listed works still see their limitation, mainly in the following three aspects: a) the constituting cells are restricted to be rectangular or cuboid, and spatial changes are enabled only along certain prescribed directions; b) smooth connections across cell boundaries cannot be guaranteed; c) the accuracy of the simulation results for the overall compliance can not be effectively predicted.

Recently, a set of novel methods were introduced, which can be roughly divided into three categories: a) the conformal-mapping-based methods \cite{vogiatzis2018computational, ye2019topology, li2020anisotropic}, b) the de-homogenisation method \cite{groen2018homogenization, groen2019homogenisation, allaire2019topology}, c) the improved asymptotic-homogenisation-based topology optimisation (AHTO plus) approach \cite{zhu2019novel, xue2020generation}.

For the conformal-mapping-based method \cite{vogiatzis2018computational, ye2019topology, li2020anisotropic}, microstructural cells are deployed, even on a manifold, guided by the use of appropriate conformal mapping functions. Consequently, structural orthotropy is preserved everywhere, provided that the constituting cells take laminate configurations. As the design functions carry geometric information now, the optimal design of the corresponding GMC becomes rather intuitive. But for an effective numerical scheme to accurately predict the compliance of generated has been barely reported until recently \cite{LiSS_arXiv2021}.

The so-called de-homogenisation method, as named by some of its initiators \cite{groen2018homogenization, groen2019homogenisation, allaire2019topology}, seeks to find the optimal material distribution with a coarse grid first. During this stage, the volume fraction of solids and the orientation of given laminate cell configurations are determined on every macroscopic pixel. The actual microstructural configuration is then resolved by projecting the chosen cells in alignment with the obtained cell orientation in a background with a much finer grid.

As for the AHTO plus method \cite{zhu2019novel, xue2020generation}, it starts with representing a GMC through a topology description function (TDF) defined by means of function composition \cite{LiuC_JAM2017, zhu2019novel}. Then asymptotic analysis, as done for asymptotic homogenisation for periodic configurations, is carried out to reformulate the original multiscale problem as a series of scale-separated problems. On the macroscale, a homogenised problem is defined as if the GMC is treated as a continuum, whose equivalent elastic moduli are computed by solving a number of cell problems defined on a fine-scale basis. Thus convergence of the derived formulation to the original fine-scale formulation is theoretically ensured. The compliance optimisation is then conducted under the framework of moving morphable components / voids (MMC/MMV) \cite{GuoX_JAM2014, ZhangWS_SMO2016}. Compared with the other two types of novel methods mentioned above, the constituting cells are also permitted to take non-rectangular/non-cuboid shape in space. Thus far, the major challenge limiting the implementation of the AHTO plus method is the computational cost, because one has to compute the microscopic cell problems as many times as the number of macroscopic finite elements, and a zoning strategy has been proposed to alleviate its generated computational burden \cite{xue2020speeding}.

The present article is thus proposed to devise a proper scheme with the use of machine learning, so as to effectively speed up the multiscale computation of the compliance of a GMC. Although the present study is carried out with a background application of the compliance design of GMC, it is actually aimed to address a ubiquitous issue, and the core question is, ``how to thoroughly accelerate the evaluations of the performance of a microstructural configuration with the simulation accuracies theoretically guaranteed at a certain level?'' A quick answer to that question, as demonstrated with the present article, is through a mutually boosting use of asymptotic analysis and machine learning. On one hand, asymptotic analysis enables one to derive scale-separation formulations converging to the original problem defined on a fine-scale basis. On doing so, the appropriate inputs and outputs to set up a machine learning model are identified, and a computational route for generating data for its training is also determined. On the other hand, machine learning enables one to accurately represent the implicit interrelationships between the homogenised quantities resulted from asymptotic analysis. The efficiency of the proposed scheme is demonstrated with (two-dimensional) numerical examples. With the homogenisation error kept below 2\%, the computational time for compliance optimisation is about 300 seconds for two-dimensional problems on a desktop computer, as the present method is used, making it as efficient as the state-of-the art projection-based method \cite{groen2018homogenization}.

In this article, we also discuss in detail a number of subtle issues that can consolidate the engagement of machine learning to asymptotic analysis. The issues involve how to ensure the positive definiteness of the homogenised elasticity tensor represented through machine learning; how to minimise the number of input arguments setting up the machine learning so as to alleviate the requirements over the size of training data; how to computationally generate data to facilitate the training of neural networks, etc. Besides, as for optimisation, we also investigate how to effectively make full use of the machine learning results to accelerate the corresponding sensitivity analysis, an issue that is highly notable for multiscale topology optimisation.

The remaining part of this article is arranged as follows. In Sec.~\ref{Sec_AHTO_intro}, the AHTO plus framework is outlined first, along with its associated challenging issues nowadays. This is followed by an introduction on using neural networks to represent the implicit function relationships identified by asymptotic analysis in Sec.~\ref{Sec_ML}, where the procedure for data generation and the training of neural networks are also discussed. After a discussion over several key issues for the numerical implementation of the present method in Sec.~\ref{Sec_preparation}, a set of numerical examples are presented to demonstrate the accuracy and efficiency of the present method in Sec.~\ref{Sec_examples}. The article concludes in Sec.~\ref{Sec_conclusion}.

\section{The AHTO plus framework\label{Sec_AHTO_intro}}
The improved asymptotic-homogenisation-based topology optimisation framework was initiated by Zhu et al. (2019) \cite{zhu2019novel} by extending the idea underlying the classical asymptotic homogenisation method for spatially periodic configurations. It comprises of three key modules: the representation of graded microstructural configuration, the calculation of the GMC compliance, and the corresponding compliance optimisation. In this section, key issues related to the AHTO plus formulation is summarised in brief.

\subsection{Topological description of graded microstructural configurations}
Given a porous design domain $\Omega$, we use $\Omega_{\text{s}}$ to denote the solid region within $\Omega$. Mathematically, $\Omega_{\text{s}}$ can be identified by introducing a topology description function $\phi(\mathbf{x})$ defined by\
\begin{equation}\label{eq2}
    \phi(\mathbf{x}) \begin{cases}\geq 0, & \mathbf{x} \in \Omega_{\text{s}}; \\ <0, & \mathbf{x} \in \Omega \backslash \Omega_{\text{s}}.\end{cases}
\end{equation}

\begin{figure}[!ht]
    \centering
    \includegraphics[width=.8\textwidth]{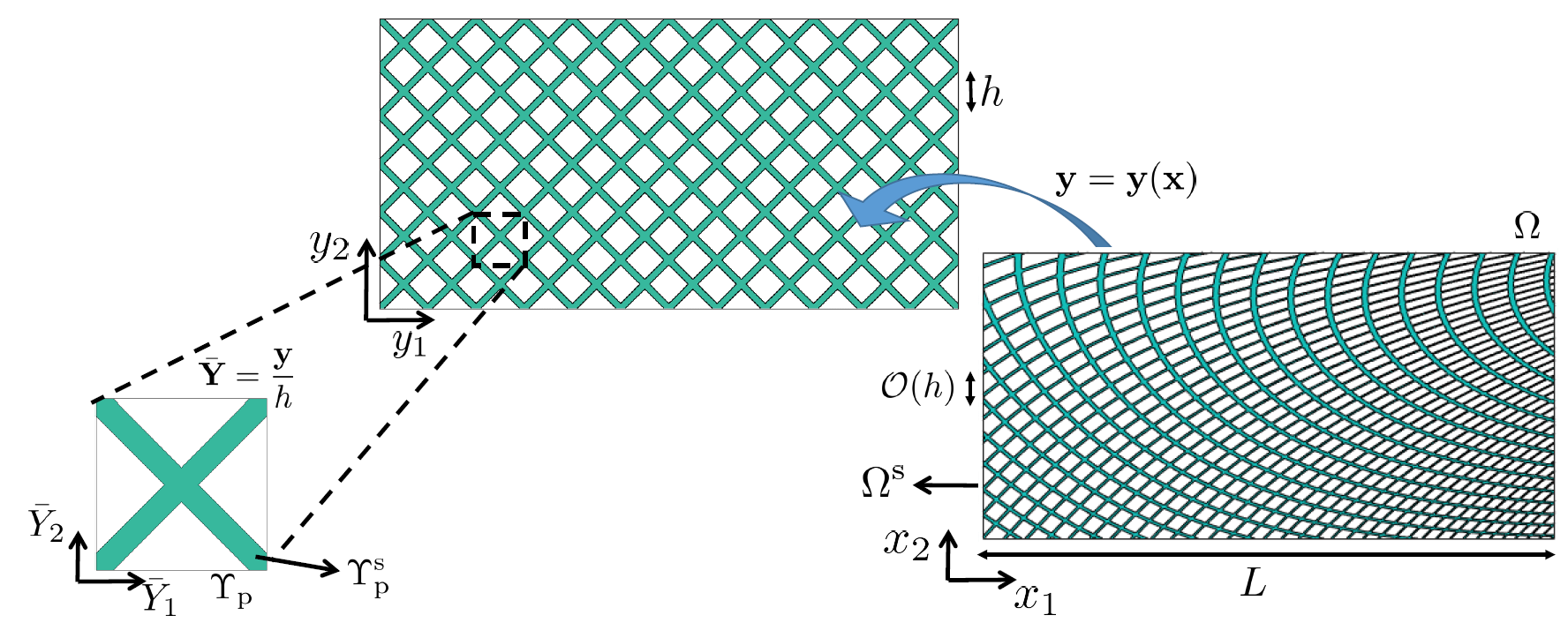}
    \caption{Generation of a graded microstructural configuration (in the panel on the right side) based on the composite topology description function given by Eq.~\eqref{eq3}. \label{f2}}
\end{figure}

Under the AHTO plus framework, a GMC is generated as follows. As shown in Fig.~\ref{f2}, a continuous mapping function $\mathbf{y}=\mathbf{y}(\mathbf{x})$ is introduced to map a GMC in an actual space (as shown in the right panel of Fig.~\ref{f2}) into a periodic configuration defined in a fictitious space (as shown in the upper panel of Fig.~\ref{f2}), which is composed of identical unit cells of size $h$.

Here we introduce a (nondimensional) rescaled coordinate system $\overline{\mathbf{Y}}$ defined by
\begin{equation}\label{eqp1}
    \overline{\mathbf{Y}}=\frac{\mathbf{y}(\mathbf{x})}{h}
\end{equation}
to measure the quantity variation in this fictitious space. Thus the unit cell, which is also termed as the ``matrix cell'' here, is defined in a nondimensional unit square/cube (as given in the bottom-left panel of Fig. \ref{f2}), i.e., $\overline{\mathbf{Y}} \in \Upsilon^{\mathrm{p}}=[-1 / 2,1 / 2]^{n}$, where $n$ represents the spatial dimension, and the superscript ``$\mathrm{p}$'' represents periodicity. Under coordinates $\overline{\mathbf{Y}}$, the microscale matrix cell, becomes nondimensional and of unit size. The structure identified by the matrix cell, with the ``X''-shape cell in Fig.~\ref{f2} as an example, can be described by a spatially periodic topology description function (defined in $\overline{\mathbf{Y}}$), as denoted by $\phi^{\mathrm{p}}(\overline{\mathbf{Y}})$. Then the TDF for the resulting GMC, as exemplified in Fig.~\ref{f2}, can be expressed by the composition of the macroscopic mapping function $\mathbf{y} = \mathbf{y}(\mathbf{x})$ and the microscopic TDF $\phi^{\mathrm{p}}(\cdot)$, i.e.,
\begin{equation}\label{eq3}
    \phi(\mathbf{x})=\phi^{\mathrm{p}}\left(\frac{\mathbf{y}(\mathbf{x})}{h}\right).
\end{equation}

The TDF given by Eq. \eqref{eq3} essentially generates a GMC by filling the design space with deformed matrix cells. The deformation operation includes stretching, rotating and twisting of the unit cells. Mathematically, such cell deformation in space is controlled by the corresponding Jacobian matrix given by
\begin{equation}\label{Jacobean}
    J_{i j}=\frac{\partial y_{i}}{\partial x_{j}}, \text{ for }i, j = 1\cdots n,
\end{equation}
which are actually formed by the (macroscopic) derivatives of the mapping function $\mathbf{y}=\mathbf{y}(\mathbf{x})$.

The advantageous features shown by the TDF~\eqref{eq3} concerning GMC generation were summarised elsewhere \cite{zhu2019novel, xue2020speeding}. It is worth noting that one can also construct a GMC from a collection of matrix cells, and the corresponding TDF, as modified from Eq.~\eqref{eq3}, reads \cite{xue2020generation}
\begin{equation} \label{TDF_multiple_cells}
\phi(\mathbf{x})=\phi^{\mathrm{p}}\left(\frac{\mathbf{y}(\mathbf{x})}{h}\right) - \zeta(\mathbf{x}),
\end{equation}
where the macroscopically defined function $\zeta(\mathbf{x})$ indicates which matrix cell is used near the (macroscopic) point $\mathbf{x}$.

\subsection{Asymptotic analysis with scales separated}
A GMC is associated with at least two length scales. In this paper, $L$ measures the size the design domain which is a representative of the macroscopic length scale, and $h$ roughly measures the size of the (deformed) matrix cells constituting the GMC which is a representative of the microscopic length scale. For a GMC, we have
\begin{equation}\label{eq1}
    \epsilon=\frac{h}{L} \ll 1.
\end{equation}
Directly analysing the mechanical properties of a GMC on the fine scale normally renders very high computational cost. It is necessary to introduce an asymptotic homogenisation method to improve the corresponding computational efficiency.

The idea is to extend the treatments of scale separation for periodic microstructure. The key difference lying in the case of GMC, is that the homogenised elastic moduli are no longer uniform in space, because of the microstructural variance in space. For greater details about the asymptotic analysis of the GMC behaviours can be found in Zhu et al. (2019) \cite{zhu2019novel}, and we simply outline the key results here.

The homogenised displacement field $\mathbf{u}^{\mathrm{H}}$ satisfies a macroscopic equilibrium equation given by
\begin{equation}\label{eq6}
    \frac{\partial}{\partial x_{j}}\left(\mathbb{C}_{i j k l}^{\mathrm{H}} \frac{\partial u_{k}^{\mathrm{H}}}{\partial x_{l}}\right)+f_{i}^{\mathrm{H}}=0, \quad \mathbf{x} \in \Omega,
\end{equation}
where $\mathbb{C}^{\mathrm{H}}=\left(\mathbb{C}_{ijkl}^{\mathrm{H}}\right)(i, j, k, l=1, \cdots, n)$ represents the elastic moduli of the equivalent continuum; $\mathbf{f}^{\mathrm{H}}=\left(f_{i}^{\mathrm{H}}\right) (i=1, \cdots, n)$ represents the homogenized body force density. Note that the Einstein summation rule is employed throughout this article. Here we let $\mathbf{f}^{\mathrm{H}}=\mathbf{0}$.

The homogenised elastic moduli $\mathbb{C}_{i j k l}^{\mathrm{H}}$ by Eq. \eqref{eq6} are calculated by
\begin{equation}\label{eq7}
    \mathbb{C}_{i j k l}^{\mathrm{H}}=\mathbb{C}_{i j k l}\left|\Upsilon_{\mathrm{p}}^{\mathrm{s}}\right|-\mathbb{C}_{i j s t} J_{n t} \int_{\Upsilon_{\mathrm{p}}^{\mathrm{s}}} \frac{\partial \xi_{s}^{k l}}{\partial \bar{Y}_{n}} \overline{\mathbf{Y}},
\end{equation}
where the third-order tensor $\xi(\overline{\mathbf{Y}}; \mathbf{x})=\left(\xi_{i}^{j k}(\overline{\mathbf{Y}} ; \mathbf{x})\right), i, j, k=1, \ldots, n$ are known as the generalised displacements, $|\Upsilon_{\mathrm{p}}^{\mathrm{s}}|$ measures the volume fraction of the solid materials within the matrix cell of interest. The third-order tensor $\xi(\overline{\mathbf{Y}}; \mathbf{x})$ is then determined by a set of cell problems (in $\overline{\mathbf{Y}}$ and parameterised with $\mathbf{x}$) governed by
\begin{equation}\label{eq8}
    J_{m j} \frac{\partial}{\bar{Y}_{m}}\left(\tilde{\mathbb{C}}_{i j k l} J_{n l} \frac{\partial \xi_{k}^{s t}}{\partial \bar{Y}_{n}}\right)=J_{m j} \frac{\partial \tilde{\mathbb{C}}_{i j s t}}{\partial \bar{Y}_{m}}, \quad \overline{\mathbf{Y}} \in \Upsilon_{\mathrm{p}},
\end{equation}
and imposed with periodic boundary conditions, where the Jacobian matrix $\mathbf{J}$ is recalled to be defined by Eq. \eqref{Jacobean}. Compared with conventional asymptotic homogenisation approaches for periodic structures, the homogenised elastic moduli $\mathbb{C}_{i j k l}^{\mathrm{H}}$ depend on the Jacobean matrix $\mathbf{J}$, which carries the information about how the matrix cell gets deformed at the macroscopic point $\mathbf{x}$.

\subsection{Compliance optimisation and sensitivity analysis \label{Sec_sensitivity}}
Upon homogenisation, the GMC compliance can be asymptotically calculated by
\begin{equation}\label{eqp3}
    \min _{\mathbf{y}(\mathbf{x}), \phi^{p}(\overline{\mathbf{Y}})} \mathcal{C}^{\text {H }}= \min _{\mathbf{y}(\mathbf{x}), \phi^{p}(\overline{\mathbf{Y}})} \int_{\Omega} \mathbb{C}_{i j k l}^{\mathrm{H}} \frac{\partial u_{i}^{\mathrm{H}}}{\partial x_{j}} \frac{\partial u_{k}^{\mathrm{H}}}{\partial x_{l}} \mathrm{~d} \mathbf{x}.
\end{equation}

For stiffness optimisation, one seeks to minimise the elastic energy of the system, and the design variables can be roughly divided into two groups: the macroscopic design variables $d_{\alpha}$ controlling the deformation of the unit cell in space and the microscopic design variables $d_{\beta}$ controlling the material topology within the unit cell.

The optimisation process can be speeded up, if the sensitives of the design variables to the elastic energy quantity $\mathcal{E}^{\mathrm{H}}$ are given. As from Xue et al. (2020) \cite{xue2020speeding}, we have
\begin{equation}\label{eqp4}
    \frac{\partial \mathcal{E}^{\mathrm{H}}}{\partial d_{\tau}}=-\int_{\Omega} \frac{\partial \mathbb{C}_{i j k l}^{\mathrm{H}}}{\partial d_{\tau}} \frac{\partial u_{i}^{\mathrm{H}}}{\partial x_{j}} \frac{\partial u_{k}^{\mathrm{H}}}{\partial x_{l}} \mathrm{~d} \mathbf{x},
\end{equation}
where $d_{\tau}$ is the collection of both the macroscopic and the microscopic design variables.

With Eq.~\eqref{eqp4}, differentiation is transformed onto the homogenised elastic moduli $\mathbb{C}_{i j k l}^{\mathrm{H}}$. Based on the composite TDF defined by Eq.~\eqref{eq3}, the macroscopic design variable $d_{\alpha}$ should be contained in the Jacobean matrix $\mathbf{J}$ defined by Eq.~\eqref{Jacobean}. Thus the chain rule reads
\begin{equation} \label{chain_rule_macro}
\frac{\partial \mathbb{C}_{i j k l}^{\mathrm{H}}}{\partial d_{\alpha}} = \frac{\partial \mathbb{C}_{i j k l}^{\mathrm{H}}}{\partial J_{pq}} \frac{\partial J_{pq}}{\partial d_{\alpha}},
\end{equation}
for $i$, $j$, $k$, $l=1$, $\cdots$, $n$. In a similar sense, the derivative of the homogenised elastic moduli with respect to the microscopic design variable $d_{\beta}$ can be calculated by
\begin{equation} \label{chain_rule_micro}
\frac{\partial \mathbb{C}_{i j k l}^{\mathrm{H}}}{\partial d_{\beta}} = \frac{\partial \mathbb{C}_{i j k l}^{\mathrm{H}}}{\partial \gamma_{\chi^{\text{p}}}} \frac{\partial \gamma_{\chi^{\text{p}}}}{\partial d_{\beta}} + \frac{\partial \mathbb{C}_{i j k l}^{\mathrm{H}}}{\partial \zeta} \frac{\partial \zeta}{\partial d_{\beta}},
\end{equation}
where $\gamma_{\chi^{\text{p}}}$ represent the parameters characterising the topology description in the unit cell of interest. Note that the second term on the right side of Eq.~\eqref{chain_rule_micro} only appears when there are more than one matrix cell with $\zeta$ recalled to be the indicator of matrix cell as defined by Eq.~\eqref{TDF_multiple_cells}.

\subsection{Major challenge\label{Ch_2_3}}
As the microstructure varies in space now, the homogenised elastic moduli become heterogeneous in the effective continuum. Consequently, one may need to solve cell problem~\eqref{eq8} as many times as the number of the finite elements for solving the homogenised problem. Facing with this challenging issue, the present article considers employing machine learning to fully release the potentials of the AHTO plus method, and the effectiveness of such a combinative use of asymptotic analysis and machine learning is then demonstrated with numerical examples in Sec.~\ref{Sec_examples}.

\section{Neural network representation of the homogenised elasticity tensor\label{Sec_ML}}
The homogenised elasticity tensor $\mathbb{C}^{\mathrm{H}}$, as calculated based on Eq.~\eqref{eq7}, is actually dependent on the third-order tensor $\xi_{i}^{j k}$ calculated from the cell problem~\eqref{eq8}. And the cell problem~\eqref{eq8} can be considered as a family of problems parameterised by the Jacobean matrix $\mathbf{J}$ and the topology within the matrix cell implied by $\phi^{\mathrm{p}}(\cdot)$. Therefore, an (implicit) function relation should be there, i.e.,
\begin{equation}\label{eqp9}
    \mathbb{C}_{i j k l}^{\mathrm{H}}=\mathbb{C}_{i j k l}^{\mathrm{H}}(\mathbf{J},\phi^{\mathrm{p}}(\cdot)),
\end{equation}
so as to relating $\mathbb{C}^{\mathrm{H}}$ to the parameters appearing in the TDF~\eqref{eq3} of a GMC.

Note that both sets of the input arguments in Eq.~\eqref{eqp9}, $\mathbf{J}$ and $\phi^{\text{p}}(\cdot)$, do not depend on the loading conditions imposed over the GMC of interest. Thus quantitative identification of the function relation~\eqref{eqp9} can be conducted at an offline stage. What makes Eq.~\eqref{eqp9} more attractive comes from the fact that $\mathbf{J}$ is defined in a macroscopic sense, and $\phi^{\text{p}}(\cdot)$ is defined within a unit cell. Hence although the evaluation of $\mathbb{C}_{i j k l}^{\mathrm{H}}$ entails considerably intensive fine-mesh calculations (also restricted in unit cells), the representation of Eq.~\eqref{eqp9} does not necessitate fine-scale meshes. This means the number of the input arguments setting up the machine learning model for use is rather limited, and the requirement on the size of the training dataset is alleviated.

In this section, we seek to use machine learning to identify this implicit, but fully determined function relationship implied by Eq.~\eqref{eqp9}.

\subsection{Identification of the inputs and outputs for machine learning}
In order to set up a machine learning model, one needs to clearly declare the corresponding inputs and outputs. Here the input/output identification is conducted following three criteria. First, the inputs and outputs should all be nondimensional and preferentially normalised. This is good for determining the ranges for data generation. Secondly, the number of inputs should be minimised, so as to keep the complexity of the trained machine learning model low. This is done by exploring the interrelationship among the input arguments. Thirdly and occasionally, certain key properties carried by the output arguments can not be preserved automatically by the trained machine learning models. Hence one may need to adjust the output arguments, so as to maintain the desired properties. An example over this issue will be presented later, where the homogenised elasticity tensor should be positively defined.

\subsubsection{Non-dimensionalisation}
Note that the expected inputs of Eq.~\eqref{eqp9}, $\mathbf{J}$ and the parameters of $\phi^{\mathrm{p}}(\cdot)$, are already nondimensionalised. Here the ranges for $\mathbf{J}$ and the parameters of $\phi^{\mathrm{p}}(\cdot)$ are quite ``normal''. This is because a too large entry of $\mathbf{J}$ may result in severe distortion of the matrix cell, which should be avoided during optimisation. Moreover, $\phi^{\mathrm{p}}(\cdot)$ is defined in a cell of unit size. Thus the parameters it carries should take values within a ``normal'' range.

Hence one just considers nondimensionalising the outputs of Eq.~\eqref{eqp9}, $\mathbb{C}^{\mathrm{H}}$, which is naturally carried out by
\begin{equation}\label{eqp11}
    \bar{\mathbb{C}}_{i j k l}^{\mathrm{H}}=\frac{\mathbb{C}_{i j k l}^{\mathrm{H}}}{E} , \text{ for }i, j, k, l = 1,\cdots, n,
\end{equation}
where $E$ is the Young's modulus of the constituting materials. Note that a bar is put atop a variable, to indicate that it is defined through nondimensionalisation.

\subsubsection{Dimension reduction for input arguments}
The inputs of Eq.~\eqref{eqp9} are from two sources: the spatially-varying deformation of the matrix cell measured by $\mathbf{J}$ and the materials structure with the matrix cell controlled by the coefficients of $\phi^{\mathrm{p}}(\cdot)$. To do dimension reduction for $\mathbf{J}$, we refer to Appendix 1, where the (nondimensional) homogenised elasticity tensor is found to be independent of $\det \mathbf{J}$, which means the determinant of $\mathbf{J}$. Hence the actual (macroscopic) inputs can be chosen to be $\mathbf{J}^{\prime}$, where one requires that
\begin{equation}\label{eqp19}
    \det \mathbf{J}^{\prime}=1.
\end{equation}

For an arbitrary $\mathbf{J}$, when calling the trained neural network, one simply rescales it by
\begin{equation}\label{eqp20}
    J_{i j}^{\prime}=(\det \mathbf{J})^{-\frac{1}{\mathrm{n}}} J_{i j}, \text{ for }i, j = 1, \cdots, n
\end{equation}
for use, where $n$ is recalled to be dimensionality number.

As for the microscopic TDF $\phi^{\mathrm{p}}(\cdot)$, there are many ways to represent it. Here the MMC framework \cite{GuoX_JAM2014} is suggested for use. The reason is twofold. First, the materials structure within a unit cell is comparably simple in general, and this favours the use of explicit topology description framework. Second, under an explicit topology description scheme, the size of (microscopic) design variables is relatively small, which keeps the required data size for training relatively low.

\subsubsection{Positive definiteness of the homogenised elasticity tensor}
The (nondimensional) homogenised elasticity tensor $\bar{\mathbb{C}}^{\mathrm{H}}$ is shown to bear symmetry as follows \cite{zhu2019novel}
\begin{equation}\label{eqp12}
    \bar{\mathbb{C}}_{i j k l}^{\mathrm{H}}=\bar{\mathbb{C}}_{j i k l}^{\mathrm{H}}=\bar{\mathbb{C}}_{k l i j}^{\mathrm{H}}, \quad \text { for } i, j, k, l=1, \cdots, n.
\end{equation}
With such symmetry, the number of outputs can be reduced from $n^4$ to $( \frac{\left(n+1\right)n}{2}+1) \frac{\left(n+1\right)n}{4}$.

Besides, the homogenised elasticity tensor is also shown positively defined, i.e. for any second-order tensor $\epsilon_{ij}$, we have
\begin{equation}\label{eqp13}
    \bar{\mathbb{C}}_{i j k l}^{\mathrm{H}} \epsilon_{ij} \epsilon_{kl} > 0,
\end{equation}
provided that $ \epsilon_{ij}  \epsilon_{ij} \neq 0$.

Note that such positive definiteness carried by $\mathbb{C}^{\mathrm{H}}$ can not be guaranteed when $\mathbb{C}^{\mathrm{H}}$ is represented by a machine learning model. To ensure such positive definiteness, we rearrange the fourth-order tensor $\mathbb{C}_{ijkl}^{\mathrm{H}}$ as a symmetric matrix $\mathbf{C}$. This is done by merging the subindex $i$ with $j$, and $k$ with $l$, i.e. $\bar{\mathbb{C}}_{ijkl}^{\mathrm{H}} = \bar{\mathbf{C}}_{IK}$, with $I,K=1,2, \cdots , \frac{n(n+1)}{2}$, and the matrix $\bar{\mathbf{C}}$ should be symmetric and positively defined, too.

Now we apply the Cholesky decomposition \cite{xu2021learning} to the matrix $\bar{\mathbf{C}}$, i.e.,
\begin{equation}\label{Cholesky}
    \bar{\mathbf{C}}=\mathbf{L L}^{\mathrm{T}},
\end{equation}
where $\mathbf{L}$ is an $\frac{n(n+1)}{2}-$ by $-\frac{n(n+1)}{2}$ lower triangular matrix, and $\mathbf{L}^{\mathrm{T}}$ denotes the transpose of $\mathbf{L}$. Then instead of $\bar{\mathbb{C}}^{\mathrm{H}}$, one may seek to determine the interrelationship between $\mathbf{L}$ and the input arguments, and the (nondimensional) homogenised elasticity tensor is calculated by Eq.~\eqref{Cholesky}.

\subsubsection{A brief summary}
To the present stage, what we expect is a surrogate model properly representing the function relationship of
\begin{equation}\label{ML_target_function}
    \mathrm{L}_{I K} = \mathrm{Z}_{I K}(\mathbf{J}^{\prime}; d_1, \cdots, d_{\beta}),
\end{equation}
where $I \geq K$ and $I, K=1, \cdots, \frac{n(n+1)}{2}$; $\det \mathbf{J}^{\prime}=1$; $d_1, \cdots, d_{\beta}$ are the parameters of the microscopic TDF $\phi^{\mathrm{p}}(\cdot)$ defined in the matrix cell $\Upsilon_{\mathrm{p}}$. As a special case, when TDF~\ref{TDF_multiple_cells} is adopted, the microscopic design variables are contracted to $\zeta$.

\subsection{Data generation}
\subsubsection{General procedure \label{Ch_3.2.1}}
Here the surrogate model of use is trained neural networks \cite{hinton1986learning}. Data are needed for training a neural network to represent Eq.~\eqref{ML_target_function}. Now the input and output arguments have been identified, and one may undertake the procedure as follows to generate the training data.

Step 1. Evaluation of the input arguments.

First, values are assigned to all the input arguments $\mathbf{J}^{\prime}$ and $d_1, \cdots, d_{\beta}$. Note that if $\mathbf{J}$ is generated arbitrarily, Eq.~\eqref{eqp20} is used to ensure that $\det \mathbf{J}^{\prime} =1$. Also note that the effective range for each input argument should be specified in advance. This issue will be further illustrated with a two-dimensional example presented later.

Step 2. Calculation of the resulting cell problem.

Evaluating $J_{i j}$ in Eq.~\eqref{eq8} by the evaluated $J^{\prime}_{i j}$ from the previous step, the governing equation for the resulting microscopic cell problem is set up. Inserting the evaluated microscopic parameters $d_1, \cdots, d_{\beta}$ into the TDF $\phi^{\mathrm{p}}(\cdot)$, the actual domain of definition for the cell problem is identified. Then with periodic boundary conditions imposed, the third-order tensor $\xi_{i}^{j k}$ can be determined base on the cell problem governed by Eq.~\eqref{eq8}.

Step 3. Calculation of the homogenised elasticity tensor.

Using Eq.~\eqref{eq7}, the (nondimensional) homogenised elasticity tensor can be evaluated.

Step 4. Application of the Cholesky decomposition.

Finally, the Cholesky decomposition~\eqref{Cholesky} is employed to calculate the lower triangular matrix $\mathbf{L}$, whose entries form the corresponding output data.

By doing so, a set of data are obtained in the multi-dimensional space spanned by the input and output arguments. The process can be repeated (in parallel) to generate more sets of data.

\subsubsection{An two-dimensional example}
A more detailed example on date generation is given with $n=2$, where the design domain becomes spatially two-dimensional.

\begin{figure}[!ht]
    \centering
    \includegraphics[width=.7\textwidth]{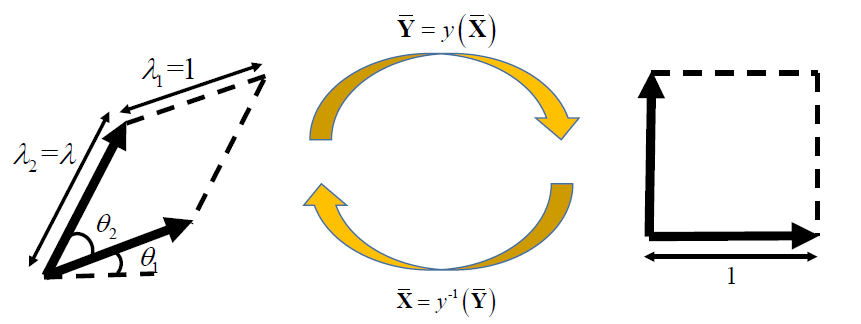}
    \caption{Based on the mapping function~\eqref{eqp1}, a parallelogram-shape cell is transformed back to the matrix cell of square shape.\label{fig_parallelogram}}
\end{figure}

Here the evaluation of $\mathbf{J}^{\prime}$ is emphasised. In theory, one may randomly assign four values to a $2-$ by $-2$ matrix, say, $\mathbf{J}$, and then calculate $\mathbf{J}^{\prime}$ with Eq.~\eqref{eqp20}. But the geometric information implied by $\mathbf{J}^{\prime}$ is fully missing then. Alternatively, we can use four other parameters to represent $\mathbf{J}$. It is noted that in two-dimensional situations, the Jacobean matrix $\mathbf{J}$ should correspond to a parallelogram configuration, which results from the deformation of the matrix cell, as shown in Fig. \ref{fig_parallelogram}.
The parallelogram can be characterised by four parameters bearing geometrical meaning: the lengths of two edges $\lambda_{1}$ and $\lambda_{2}$, a rational angle $\theta_{1}$ and an interior angle $\theta_{2}$. Since similar expansion of the parallelogram does not give rise to changes in the corresponding homogenised elasticity tensor, as demonstrated in Appendix 1, we simply let $\lambda_{1}=1$, and $\lambda_{2}=\lambda$.

With these geometrically meaningful parameters, the Jacobean matrix $\mathbf{J}$ identifying the parallelogram in Fig.~\ref{fig_parallelogram} is formulated by \cite{zhu2019novel}
\begin{equation}\label{eq24}
    \mathbf{J}=\frac{1}{\lambda \sin \theta_{2}}\left(\begin{array}{cc}
        \lambda \sin \left(\theta_{1}+\theta_{2}\right) & -\lambda \cos \left(\theta_{1}+\theta_{2}\right) \\
        -\sin \theta_{1} & \cos \theta_{1}
    \end{array}\right).
\end{equation}

Therefore, the evaluation of $\mathbf{J}$ can be done as follows. We first randomly assign certain values to the three geometrical parameters $\lambda, \theta_{1}$ and $\theta_{2}$. Then we can evaluate $\mathbf{J}$ through Eq.~\eqref{eq24}, and the values of $\mathbf{J}^{\prime}$ are finally obtained with the use of Eq.~\eqref{eqp20}.

Such a way of evaluating $\mathbf{J}^{\prime}$ is good for keeping their values within a reasonable range. As discussed above, highly distorted cells are not expected for accurate design. Geometric representation of a deformed matrix cell offers a way to explicitly control its shape. Here we simply let
\begin{equation}\label{eq25}
    \frac{1}{\lambda_{\max }} \leq \lambda \leq \lambda_{\max }, \quad 0 \leq \theta_{1}<2 \pi, \quad \theta_{\min } \leq \theta_{2} \leq \pi-\theta_{\min },
\end{equation}
when $\mathbf{J}$ get evaluated.

In order to ensure an even distribution of the obtained data points in the space spanned by the input arguments domain, the concept of low discrepancy sequence \cite{mishra2021enhancing} is employed for data generation. Here the so-called Sobol sequence is used to generate random data sequence in the three-dimensional space spanned by $\lambda$, and $\theta_{1}$ and $\theta_{2}$. The actual values of the input argument $\mathbf{J}'$ are then obtained following Eqs.~\eqref{eq24} and \eqref{eqp20}.

As for data generation, the microscopic inputs $d_1$, $\cdots$, $d_{\beta}$ in Eq.~\eqref{ML_target_function} should be activated, as further discussed in Sec.~\eqref{Sec_uniform_load}. But here for illustration, we fix the matrix cell to be ``X''-shape, as shown in Fig.~\ref{f2}, bearing a volume fraction of 30\%. Thus the microscopic parameters $d_1$, $\cdots$, $d_{\beta}$ in Eq.~\eqref{ML_target_function} are fixed (temporarily), and the total number of active input variables for the function relation~\eqref{ML_target_function} is reduced to three, that is, the three free components of $\mathbf{J}^{\prime}$. Now following the procedure detailed in Sec.~\ref{Ch_3.2.1}, 4,000 data points are generated. Note that although 4,000 data points do not sound many in a context of machine learning, it is shown by the our numerical examples in Sec.~\ref{Sec_examples} that the trained neural networks can deliver satisfactory accuracy underpinning the successive two-dimensional stiffness optimisation (with only the macroscopic design variables activated). The reason behind such a low requirement over the size of training data is perhaps due to the generic smoothness carried by the function relation~\eqref{ML_target_function}. As a result, the hidden trends against the input variables should be captured with the use of a limited amount of data.

\subsection{Machine learning}

\subsubsection{Settings and general performances \label{Ch_3.3.1}}
The back-propagation (BP) neural network is adopted for representing the relationship \eqref{ML_target_function}, and the networks employed here are set up in a fully connected manner. More details about the setting for neural network training are compiled in Table~\ref{Table_ML_setting}.
\begin{table}[!ht]
    \centering
    \renewcommand{\tablename}{Table}
    \caption{General settings for neural network training.}
    \setlength{\tabcolsep}{7mm}{\begin{tabular}{cc}
            \toprule
            Parameter & Value/Method \\
            \midrule
            Optimiser & Levenberg-Marquardt \\
            Activation function & Sigmoid \\
            Maximum number of training & 1000 \\
            Learning rate & $0.01$ \\
      %      Target mean square error & 0 \\		
            \bottomrule
    \end{tabular}}
    \label{Table_ML_setting}	
\end{table}

For two-dimensional cases, there are six independent outputs corresponding to the six entries of the lower triangular matrix $\mathbf{L}$ defined by Eq.~\eqref{Cholesky}. They may be represented by a single neural network or in separate. It is found that describing them by six independent networks bearing simple structures seems to deliver sufficient accuracy at a low training cost. Hence six individual neural networks are trained with their structural profiles summarised in Table~\ref{tbl:table-1}.
\begin{table}[!ht]
    \centering
    \renewcommand{\tablename}{Table}
    \caption{The structures of the neural networks trained for representing the six entries of the lower triangular matrix $\mathbf{L}$ defined by Eq.~\eqref{Cholesky}. All the neural networks bear no more than 3 hidden layers, and the number of neurons in each layer are given.}
\begin{tabular}{|l|ccc|}
    \hline
     \diagbox{Net}{Number}{Layer} & Layer~1 & Layer~2 & Layer~3 \\
    \hline
    $L_{11}$ & 32 & 24 & 9\\
    $L_{21}$ & 32 & 24 & 9\\
    $L_{31}$ & 30 & 30 & - \\
    $L_{22}$ & 20 & 20 & - \\
    $L_{32}$ & 32 & 24 & 9\\
    $L_{33}$ & 30 & 20 & 15\\
    \hline	
\end{tabular}
    \label{tbl:table-1}
\end{table}

\begin{figure}[!ht]
    \centering
    \subfigure[$L_{11}$]{\includegraphics[width=.28\textwidth]{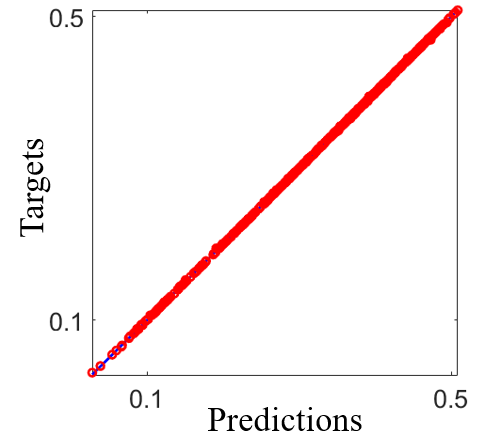}}
    \subfigure[$L_{21}$]{\includegraphics[width=.28\textwidth]{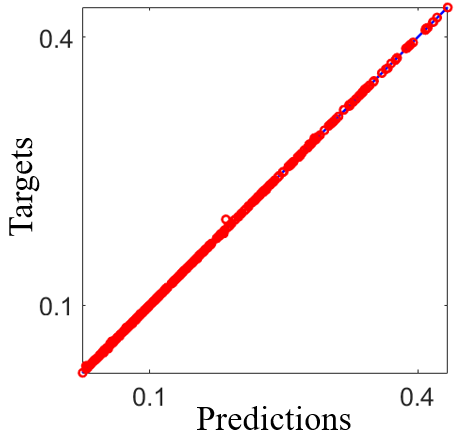}}
    \subfigure[$L_{31}$]{\includegraphics[width=.28\textwidth]{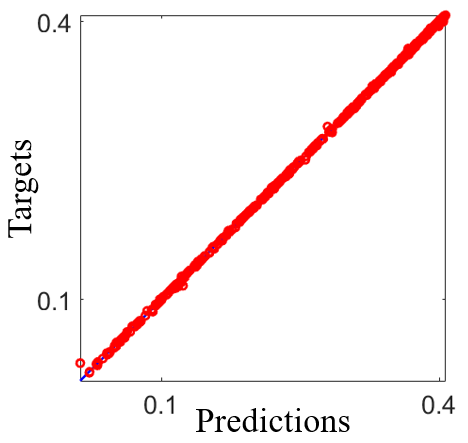}}
    \subfigure[$L_{22}$]{\includegraphics[width=.28\textwidth]{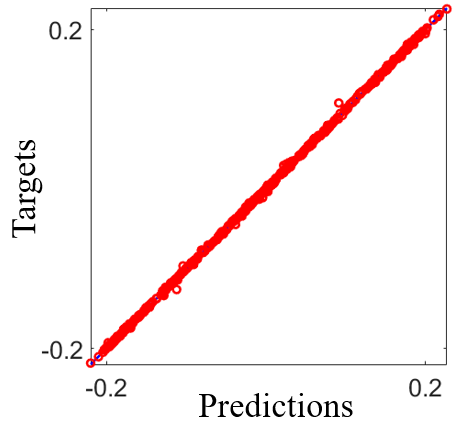}}
    \subfigure[$L_{32}$]{\includegraphics[width=.28\textwidth]{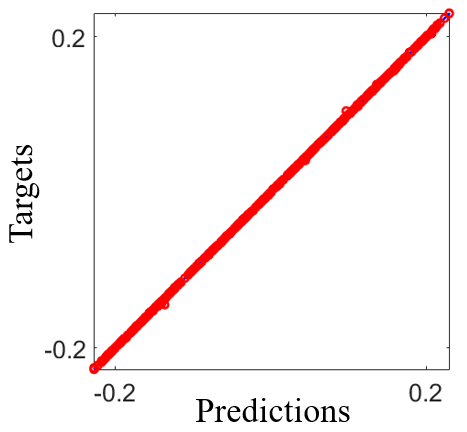}}
    \subfigure[$L_{33}$]{\includegraphics[width=.28\textwidth]{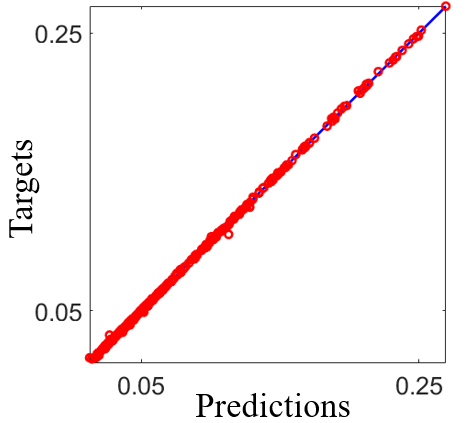}}
    \caption{General performances of the trained neural networks for the six entries of the lower triangular matrix $\mathbf{L}$ defined by Eq.~\eqref{Cholesky}\label{f7}}
\end{figure}

Among the 4,000 sets of data, 500 sets are randomly picked up to form a verification set, against which the root-mean-square error (RMSE) values no greater than 0.4\% are recorded for all the six neural networks. In Fig.~\ref{f7} shows how the predictions are drawn against their targets for all entries of $\mathbf{L}$.
Almost all the data points are found sitting on the line of $y=x$, indicating almost perfect predictions from the trained neural networks.

Several reasons may give rise to such excellent training performances with simply 4,000 sets of data. First, the function relation~\eqref{ML_target_function} bears sufficient smoothness with its input arguments. Second, the number of input arguments here is 3, which is rather low if viewed in a machine learning perspective. However, it may also be demonstrated here that proper quantitative analysis over problems even presumed to be represented with machine learning models is of great value, in the sense that certain hidden functional relations can be explored and the dimensionality of the data space is somehow minimised.

\subsubsection{Prediction of rotational symmetry}
It is noted that the homogenised elasticity tensor $\mathbb{C}^{\mathrm{H}}$ bears rotational symmetry in theory. To be precise, if we only change the value of $\theta_1$ as in Fig.~\ref{fig_parallelogram}, the values of $\mathbb{C}^{\mathrm{H}}$ are expected to be rotationally invariant. Nevertheless, a neural network is not ``aware'' of the existence of such symmetry. Therefore, as a means for testing the performance of the trained neural network, we examine whether such rotational symmetry can be reproduced by the trained neural network.

\begin{figure}[!ht]
  \centering
  \subfigure[RMSE=0.062\%]{\includegraphics[width=.32\textwidth]{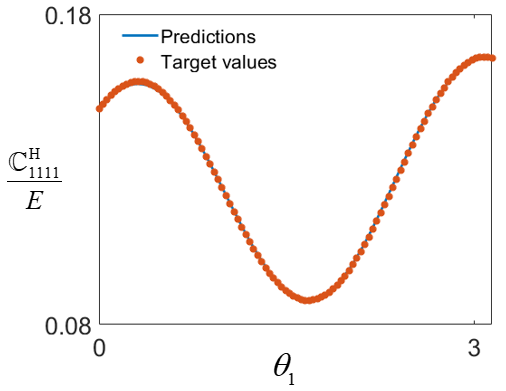}}
  \subfigure[RMSE=0.034\%]{\includegraphics[width=.32\textwidth]{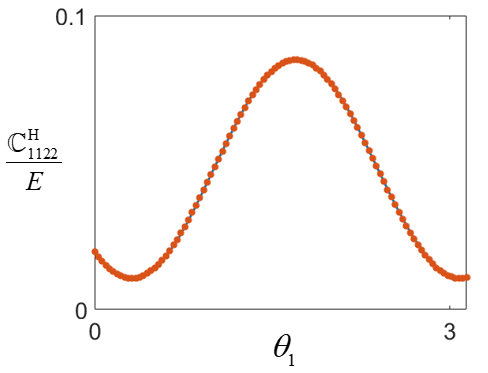}}
  \subfigure[RMSE=0.18\%]{\includegraphics[width=.32\textwidth]{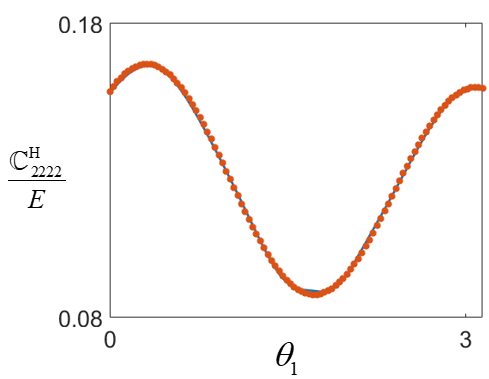}}
  \subfigure[RMSE=0.12\%]{\includegraphics[width=.32\textwidth]{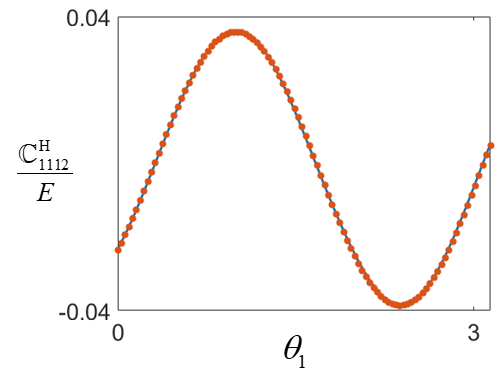}}
  \subfigure[RMSE=0.15\%]{\includegraphics[width=.32\textwidth]{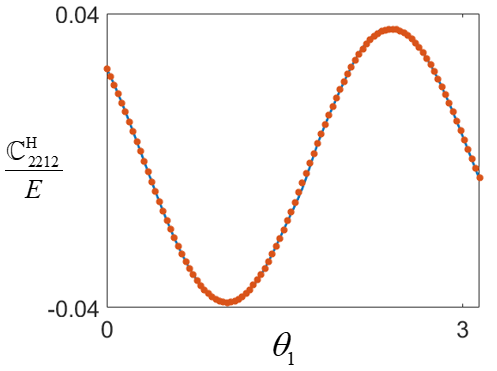}}
  \subfigure[RMSE=0.15\%]{\includegraphics[width=.32\textwidth]{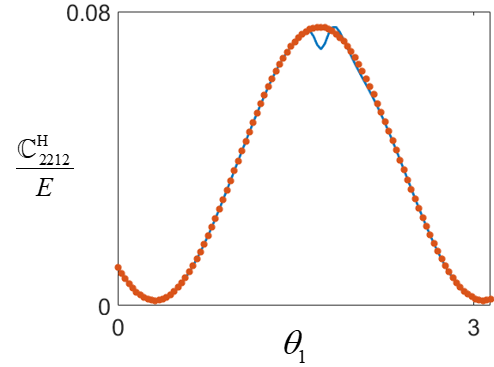}}
    \caption{Examination over the rotational symmetry supposed to be carried by the homogenised elasticity tensor predicted by the trained neural networks. \label{f6}}
\end{figure}

To this end, we fix $\lambda=1$ and $\theta_{2}=\pi/2$, and assign values to $\theta_{1}$ in an equidistant manner between 0 and $\pi$. Then we evaluate a set of $\mathbb{C}_{ijkl}^{\mathrm{H}}$ with $\theta_1=0$ following the procedure detailed in Sec.~\ref{Ch_3.2.1}, and evaluate $\mathbb{C}_{ijkl}^{\mathrm{H}}$ in the other cases based on the mentioned rotational symmetry. These 100 obtained data are used as the target data, which are compared with the corresponding results predicted by the trained neural networks. The comparative results are summarised by Fig. \ref{f6}.
It is read that the RMSE values are no greater than $0.2\%$, demonstrating that the supposed rotational symmetry of $\mathbb{C}^{\mathrm{H}}$ has been effectively ``learned'' by the neural network after training.

Accurate predictions over $\mathbb{C}^{\mathrm{H}}$ should also benefit the successive optimisation exercise. This is because the differentiation calculations~\eqref{chain_rule_macro} and \eqref{chain_rule_micro} underpinning sensitivity analysis can be conducted with a central difference scheme with confidence. This issue will be revisited in Sec.~\ref{Sec_sensitivity_CF}.

\section{Stiffness optimisation of GMC\label{Sec_preparation}}
The formulation with a combinative use of asymptotic analysis and machine learning has been derived. In this section, we discuss the critical issues on implementing it for the stiffness optimisation of GMC.

\subsection{Design variables\label{Ch_4.1}}
The design variables used for representing a GMC here can be divided into two categories. One group is formed by the controlling parameters of the macroscopic mapping function $\mathbf{y}(\mathbf{x})$, which casts its influence through the Jacobian matrix $\mathbf{J}$. The other group is formed by the microscopic design variables, which come from the microscopic TDF $\phi^{\text{p}}(\cdot)$ defined within the matrix cell $\Upsilon_{\text{p}}$. In this article, we are mainly focused on controlling the macroscopic gradient features of a GMC. Thus the matrix cell here is set fixed to be the ``X''-shape configuration as shown in the lower left panel of Fig.~\ref{f2}. Hence only the macroscopic design variables are activated here.

It is worth pointing out that the macroscopic design variables may be expressed in two ways. One is made in a geometric manner, where $\theta_{1}, \theta_{2}$ and $\lambda$ shown in Fig.~\ref{fig_parallelogram} are adopted for optimisation. The other is through an algebraic representation, where the mapping function $\mathbf{y}(\mathbf{x})$ is expressed directly, and the corresponding design variables are the parameters controlling it.  Although it is more intuitive to describe the design variables in a geometric way, challenges exist due to the following reasons.

From Eq.~\eqref{eq24}, we can obtain the values of $\mathbf{J}$ given $\theta_{1}, \theta_{2}$ and $\lambda$. Since $\mathbf{J}$ is derived by taking the partial derivatives of $\mathbf{y}$, it is not easy to calculate $\mathbf{y}$ by integrating $\mathbf{J}$, as the corresponding total differentiation condition must be satisfied. However, the inter-relation among $\theta_{1}, \theta_{2}$ and $\lambda$ is generally not analytic, except for special circumstances (for example, in two-dimensional cases where $\mathbf{y}(\mathbf{x}) $ is confined to be a conformal mapping \cite{allaire2019topology, LiSS_arXiv2021}). Thus imposing the integrability condition is not a trial task, and this poses great difficulties in calculating $\mathbf{y}(\mathbf{x})$, which is needed for generating the resulting GMC. Therefore, although not intuitive enough, the parameters of $\mathbf{y}(\mathbf{x})$ are adopted as the macroscopic design variables in this article.

The mapping function $\mathbf{y}(\mathbf{x})$ can be parameterised in the following form:
\begin{equation}\label{eq43}
    \mathbf{y}(\mathbf{x})=\sum_{\tau=1}^{K} h_{\tau} \phi_{\tau}(\mathbf{x}),
\end{equation}
where $\phi_{\tau}(\mathbf{x})$ represents a series of basis functions, $\tau=1, \cdots, K$; $h_{\tau}$ is the corresponding coefficient. Therefore, the mapping function can be expressed in a fully explicit form. And the spatial change is now completely controlled by the values of $h_{\tau}$. Here the polynomial mapping function is adopted, i.e.,
\begin{equation}\label{eq44}
    y_{i}=a_{i j} x_{j}+\frac{1}{2} b_{i j k} x_{j} x_{k}+\frac{1}{3} c_{i j k l} x_{j} x_{k} x_{l},
\end{equation}
$i=1, \cdots, n$ with the following symmetry conditions to hold: $b_{i j k}=b_{i k j}, c_{i j k l}=c_{ijl k}=c_{i lk g}=c_{i k l}$.

With Eq.~\eqref{Jacobean}, the Jacobian matrix $\mathbf{J}$ is then expressed by
\begin{equation}\label{J_expression}
    J_{i j}(\mathbf{x})=a_{i j}+b_{i j k} x_{k}+c_{i j k l} x_{k}x_{l},
\end{equation}
$i,j=1, \cdots, n$. Therefore, the design variables obtained for the optimisation of the macroscopic problem are $a_{i j}$, $b_{i j k}$ and $c_{i j k l}$.

\subsection{Constraints}
The asymptotic results for $\mathbb{C}^{\mathrm{H}}$ become inaccurate, when the matrix cell is highly distorted in space. As suggested by Eq.~\eqref{eq25}, $\lambda$ and $\theta_2$ shown in Fig.~\ref{f6} should be kept within a medium range. When the training data are generated for the examples shown here, we require $1/3<|\lambda|<3$ and $\pi/4<\theta_2<\pi/3$. These are actually the effective ranges for appropriate use of the trained neural networks. During optimisation, constraints should be imposed to ensure that these conditions should be met everywhere in the GMC. It is suggested from the numerical examples later that the optimisation is more susceptible to the limits on $\theta_2$. Here we impose a more strict requirement, $\theta_{2} \in (\frac{\pi}{4}, \frac{3 \pi}{4} )$. Mathematically, the constraints read
\begin{equation}\label{constraints0}
    \frac{1}{9} \leq \lambda^{2} \leq 9, \quad \sin \left(\theta_{2}\right) \geq \frac{\sqrt{2}}{2}.
\end{equation}

As mentioned in Sec.~\ref{Ch_4.1}, the actual design variables are the parameters of $\mathbf{y}(\mathbf{x})$ or the entries of $\mathbf{J}$. With Eq.~\eqref{eq24}, they are related to the geometric parameters appearing in Eq.~\eqref{constraints0} by
\begin{equation}\label{eq34}
    \lambda^{2}=\frac{J_{11}^{2}+J_{12}^{2}}{J_{21}^{2}+J_{22}^{2}}, \quad \sin \left(\theta_{2}\right)=\frac{J_{11} J_{22}-J_{12} J_{21}}{\sqrt{J_{11}^{2}+J_{12}^{2}} \sqrt{J_{21}^{2}+J_{22}^{2}}}.
\end{equation}
This means that the constraints given by Eq.~\eqref{constraints0} turn implicit against the design variables. Hence the $p-\text{norm}$ constraint method needs to be adopted. To this end, we introduce
\begin{subequations} \label{f_Heviside}	
    \begin{equation}\label{eq32}
        f_{1}=\mathrm{H}\left[\left[\int_{\Omega}\left(\lambda^{2 p}+\frac{1}{\lambda^{2 p}}\right) \mathrm{d} \mathbf{x}\right]^{\frac{1}{p}}-9\right],
    \end{equation}
and
    \begin{equation}\label{eq33}
        f_{2}=\mathrm{H}\left[\left[\int_{\Omega}\left(\frac{1}{\sin \theta_{2}}\right)^{p} \mathrm{d} \mathbf{x}\right]^{\frac{1}{p}}-\sqrt{2}\right],
    \end{equation}	
\end{subequations}
where $\mathrm{H}(\cdot)$ is the Heaviside function satisfying
\begin{equation}\label{eqp16}
    \mathrm{H}(t)= \begin{cases}1, & t>0 ; \\ 0, &  \text { otherwise }.\end{cases}
\end{equation}

Here $p$ is should be large, because
\begin{equation}\label{eqp17}
    \lim_{p \rightarrow \infty}\left[\int_{\Omega}\left(\lambda^{2 p}+\frac{1}{\lambda^{2 p}}\right) \mathrm{d} \mathbf{x}\right]^{\frac{1}{p}}=\max _{\mathbf{x} \in R}\left[\max \left(\lambda^{2}, \frac{1}{\lambda^{2}}\right)\right].
\end{equation}
Then $f_1$ and $f_2$ defined by Eqs.~\eqref{f_Heviside} are expected to vanish, when constraints~\eqref{constraints0} are met. Hence the cost for violating the constraint should be quite high, if $f_1$ and $f_2$ both multiplied by large numbers are added to the target function for minimisation.

\subsection{Optimisation formulation}
In summary, the optimisation problem is established by
\begin{equation}\label{eq31}
    \begin{aligned}
        &\text { Find } \quad \mathbf{y}=\mathbf{y}(\mathbf{x}) \in \mathcal{U}_{\mathrm{y}}, \phi(\overline{\mathbf{Y}}) \in \mathcal{U}\left(\mathrm{Y}_{\mathrm{p}}\right) \\
        &\text { Minimise } \quad \mathcal{F}=\int_{\Omega^{I}} \mathbb{C}_{i j k l}^{\mathrm{H}} \frac{\partial u_{i}^{\mathrm{H}}}{\partial x_{j}} \frac{\partial u_{k}^{\mathrm{H}}}{\partial x_{l}} \mathrm{~d} \mathbf{x}+K_{1} f_{1}+K_{2} f_{2} \\
        &\int_{\Omega} \mathbb{C}_{i j k l}^{\mathrm{H}} u_{i, j} v_{k, l} d \mathbf{x}=\int_{\Omega} f_{i}^{\mathrm{H}} v_{i} \mathrm{~d} \mathbf{x}+\int_{\Gamma_{\mathrm{t}}} t_{i} v_{i} \mathrm{~d} S \\
        &\forall \mathbf{v} \in \mathcal{U}_{a d} \\
        &\mathbb{C}_{i j k l}^{\mathrm{H}}=E \mathbf{L} \mathbf{L}^{ \mathrm{~T}} \quad, I=1, \cdots, N \\
        &\mathbf{L}=\mathbf{Z}\left(\mathbf{J}^{\prime}(\mathbf{x})\right) \quad, I=1, \cdots, N ; \\
        &\mathbf{u}^{\mathrm{H}}=\overline{\mathbf{u}}, \quad \text { on } \Gamma_{\mathrm{u}} \\
        &V_{\mathrm{f}} \leq \bar{V}
    \end{aligned}
\end{equation}
where $\mathbf{v}$ is the virtual displacement field; $\mathcal{U}_{\mathrm{y}}$ and $\mathcal{U}\left(\mathrm{Y}_{\mathrm{p} }\right)$ are the macro mapping function space and the micro topological description function space, respectively; $K_{1}$ and $K_{2}$ are two large numbers; $f_1$ and $f_2$ are given by Eqs.~\eqref{constraints0}; $\mathbf{Z}(\cdot)$ is the function represented by the trained neural network; $\mathbf{L}$ is the lower triangular matrix obtained by Cholesky decomposition.

\subsection{Sensitivity analysis\label{Sec_sensitivity_CF}}
The machine learning representation of the homogenised elasticity tensor $\mathbb{C}^{\mathrm{H}}$ also benefits the numerical implementation of the sensitivity analysis results given in Sec.~\ref{Sec_sensitivity}.

As from the sensitivity formulation~\eqref{chain_rule_macro} for the macroscopic design variables, for example, we want to calculate
\begin{equation}\label{chain_rule_later}
    \frac{\partial \mathbb{C}_{i j k l}^{\mathrm{H}}}{\partial d_{\alpha}}=\frac{\partial \mathbb{C}_{i j k l}^{\mathrm{H}}}{\partial J_{p q}^{\prime}} \frac{\partial J_{p q}^{\prime}}{\partial d_{\alpha}},
\end{equation}
for $i$, $j$, $k$, $l$, $p$ and $q=1$, $\cdots$, $n$, for the numerical examples presented here. As the parameterisation of $\mathbf{J}'$ is explicit (when Eq.~\eqref{J_expression} is used), the calculation of $\frac{\partial J_{p q}^{\prime}}{\partial d_{\alpha}}$ is straightforward, and the time-demanding part should come from the evaluations of $\frac{\partial \mathbb{C}_{i j k l}^{\mathrm{H}}}{\partial J_{p q}^{\prime}}$ in Eq.~\eqref{chain_rule_later}.

In existing asymptotic homogenisation frameworks, this is done by employing the adjoint method \cite{liu2002mapping}, where $\frac{\partial \mathbb{C}_{i j k l}^{\mathrm{H}}}{\partial J_{p q}^{\prime}}$ are calculated by means of integrals over
the matrix cell of interest \cite{xue2020speeding}. The computation of these integrals, albeit being affordable, is actually quite time-consuming. But now, with $\mathbb{C}^{\mathrm{H}}$ represented through neural network, we can simply implement a central difference scheme to calculate
\begin{equation}\label{central_diff}
    \frac{\partial \mathbb{C}_{i j k l}^{\mathrm{H}}}{\partial J_{p q}^{\prime}} \approx \frac{\mathbb{C}_{i j k t}^{\mathrm{H}}\left(\left.\mathbf{J}^{\prime}\right|_{J_{p q}^{\prime}+\Delta J_{p q}}\right)-\mathbb{C}_{i j k}^{\mathrm{H}}\left(\left.\mathbf{J}^{\prime}\right|_{J_{p q}^{\prime}-\Delta J_{p q}}\right)}{2 \Delta J_{p q}}
\end{equation}
for $i$, $j$, $k$, $l$, $p$, $q=1$, $\cdots$, $n$. Given the demonstrated performance of the trained neural network in Sec.~\ref{Sec_ML}, the approximating accuracy with the use of Eq.~\eqref{central_diff} should be sufficiently high. With the sensitivity quantities calculated, the algorithm of moving asymptotes (MMA) \cite{Svanberg_IJNME1987} is employed as the optimiser. The effectiveness of the sensitivity analysis conducted in this way is also addressed with the numerical examples shown in Sec.~\ref{Sec_examples}.

\subsection{Issues for speeding the proposed algorithm up\label{Ch_3.4}}
\begin{figure}[!ht]
    \centering
    \includegraphics[width=.6\textwidth]{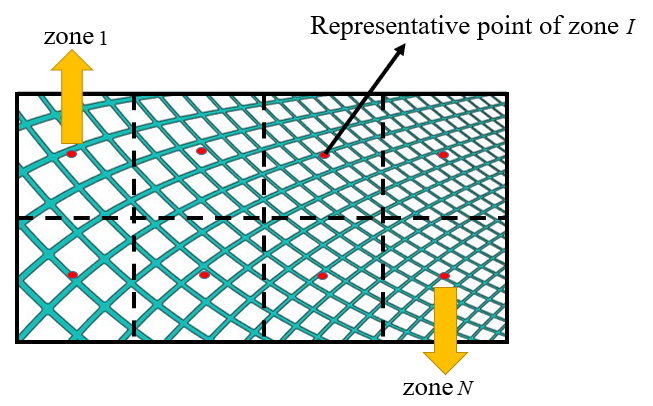}
    \caption{The concept of the zoning strategy \cite{xue2020speeding}.\label{f8}}
\end{figure}

To further speed up the computation with the homogenisation accuracy maintained at a satisfactory level, the idea of a zoning strategy is adopted \cite{xue2020speeding}, where the entire design domain can be divided into $N$ subdomains, as shown in Fig.~\ref{f8}.
One simply uses the homogenised elasticity tensor evaluated at a representative point to represent that for the whole subdomain. Such representative points are normally chosen as the central points of the subdomains. Apparently, the more the whole domain is partitioned, the more accurate the computational results are. It has been shown by Xue et al. (2020) \cite{xue2020speeding} that for two-dimensional stiffness optimisation, when each subdomain contains no more than 25 (macroscopic) finite elements, the accuracy of the compliance calculation is high enough to be used for optimisation. The numerical examples presented in Sec.~\ref{Sec_examples} also demonstrate the stableness of using the zoning strategy.

The use of such a zoning treatment also enables further acceleration as the calling of the trained neural network can be run in parallel. This is because the evaluation for the homogenised elasticity tensor $\mathbb{C}^{\mathrm{H}}$ in one subdomain is fully independent on that in the other subdomains.

\begin{figure}[!ht]
    \centering
    \includegraphics[width=.8\textwidth]{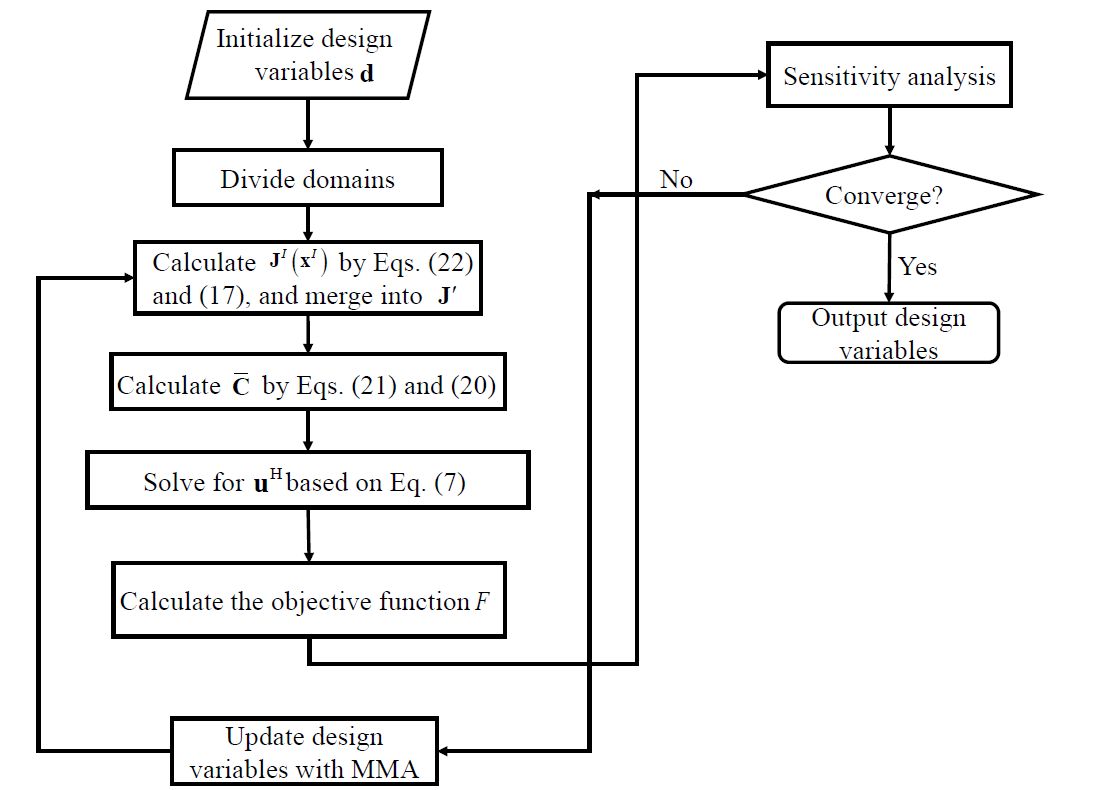}
    \caption{The flowchart for carrying out the stiffness optimisation algorithm proposed in the article.\label{f9}}
\end{figure}

In summary, the flowchart of present algorithm numerical implementation is shown in Fig.~\ref{f9}.

\section{Numerical examples and discussion\label{Sec_examples}}
In this section, we demonstrate both accuracies and efficiencies enabled by using the present homogenisation formulation supported by machine learning. To this end, a number of representative (two-dimensional) examples are presented. To address the solution accuracies, the homogenised results are compared with the underlying fine-scale results. To address the solution efficiencies, the optimisation (for two-dimensional cases) is found to cost roughly 400 seconds in total on a standard desktop computer with a singly activated core. The advantages exhibited by the present method are also shown in comparison with the original zoning strategy \cite{xue2020speeding}. Here the overall design domain is occupied by a cantilever short mean with its right edge fixed as shown in Fig.~\ref{f10}. The domain is of size $2\times1$.

\subsection{The case under uniform loading\label{Sec_uniform_load}}

\begin{figure}[!ht]
    \centering
    \includegraphics[width=.5\textwidth]{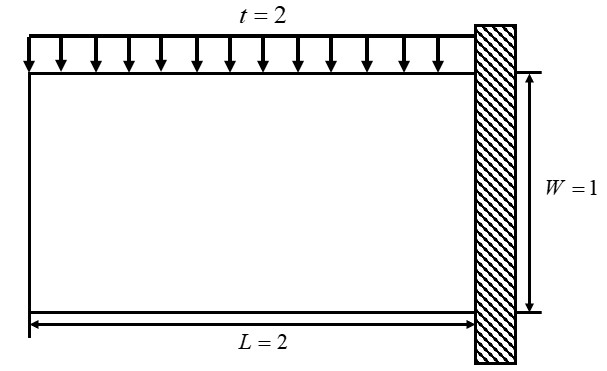}
    \caption{The case with a uniformly distributed load applied on top of a cantilever short beam.\label{f10}}
\end{figure}

First we consider the case, as shown in Fig.~\ref{f10}, where the short beam withstands a uniform load of magnitude 2 applied on its top. We seek to fill the design domain by distorting the ``X''-shape matrix cells as shown in left panel of Fig.~\ref{f2}, and the volume fraction is set to be 30\%.

\begin{figure}[!ht]
    \centering
    \subfigure[]{\includegraphics[width=.64\textwidth]{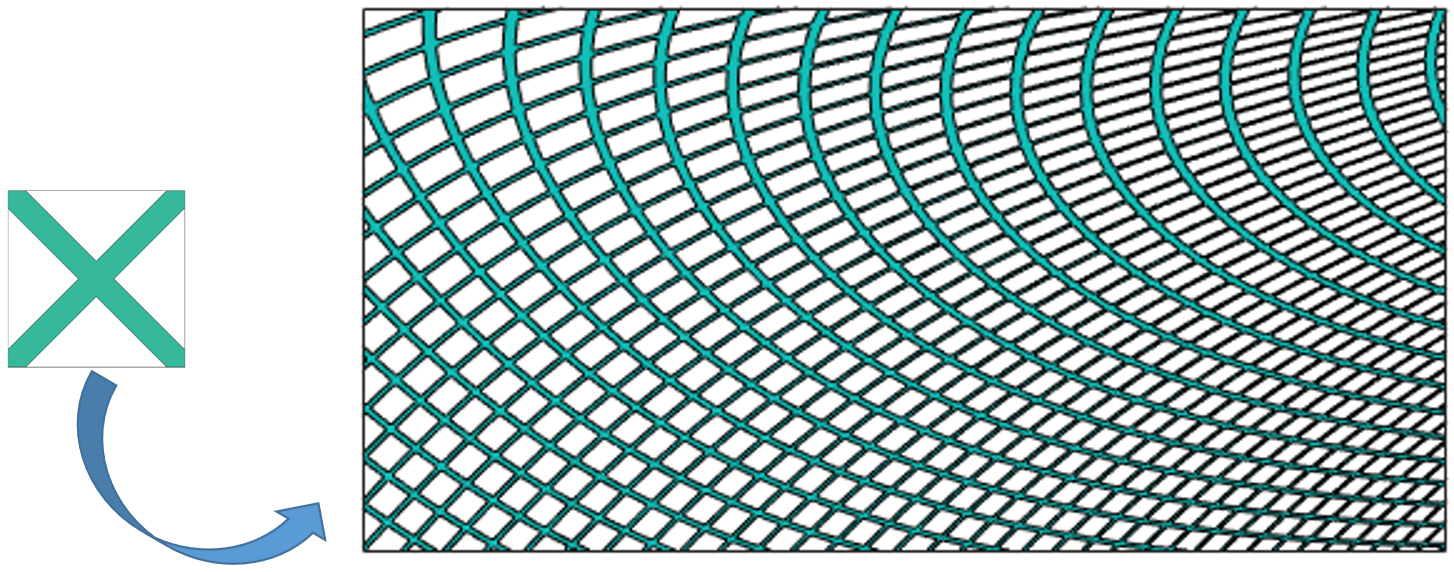}}
    \subfigure[]{\includegraphics[width=.32\textwidth]{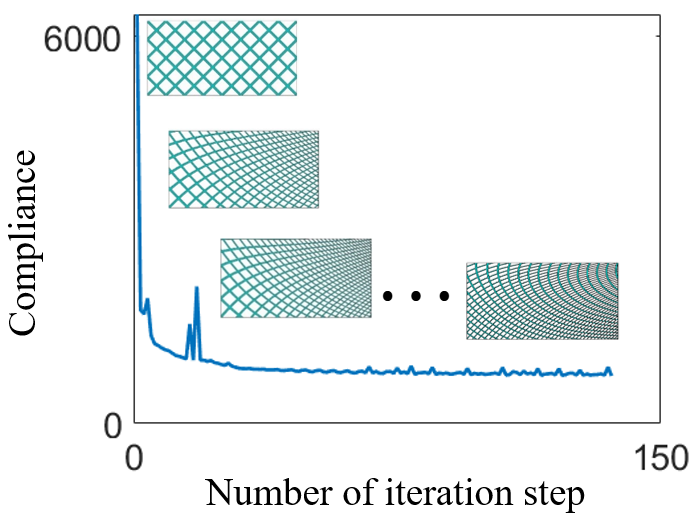}}
    \caption{Optimised results for the uniform loading case shown in Fig.~\ref{f10}. (a) The optimised GMC configuration with the matrix cell marked on the left side. (b) The evolution of the GMC compliance during optimisation. \label{fig_optimised_uniform_load}}
\end{figure}

The optimised GMC based on the present scheme is shown by Fig.~\ref{fig_optimised_uniform_load}(a), with the optimisation procedure summarised in Fig.~\ref{fig_optimised_uniform_load}(b).
To ensure the convergence of the homogenised results to the actual situations, the optimised configuration is also analysed with an extremely fine grid ($1600\times800$), where the microstructural details can be resolved. Note that to facilitate the fine-mesh computation, a thin solid layer of thickness 1/200 is attached to the top of the design domain where the beam is subject to loads. In this scenario, the homogenisation computation predicts the GMC compliance to be 710.6, which is just 0.83\% deviation from the fine-mesh calculation which predicts the compliance to be 704.8. This effectively demonstrates the accuracy of the present homogenisation approach supported by the trained neural networks.

The same problem has also been considered using the original zoning strategy \cite{xue2020speeding}, where the optimised GMC then bears a compliance of 852.5. This means that the stiffness of the optimised GMC is improved by 13\% using the present scheme equipped with machine learning models. Such improvement may be due to the fact that the compliance calculation is presumably more accurate when the homogenised elasticity tensor is represented by machine learning models. It is shown in detail in Sec.~\ref{Sec_time} that the use of machine learning enables one to do computation (on a desktop computer) with 5,000 subdomains, which is far greater than what can be done with the original zoning strategy (at most 512 subdomains on a computer containing 16 processors). An increase in the number of subdomains boosts effective calculation for the homogenised stress field in the GMC, and the evaluation of the compliance is thus far more accurate. Greater details about the efficiency of the present method due to the use of machine learning will be given in Sec.~\ref{Sec_time}.

\begin{figure}[!ht]
    \centering
    \includegraphics[width=.7\textwidth]{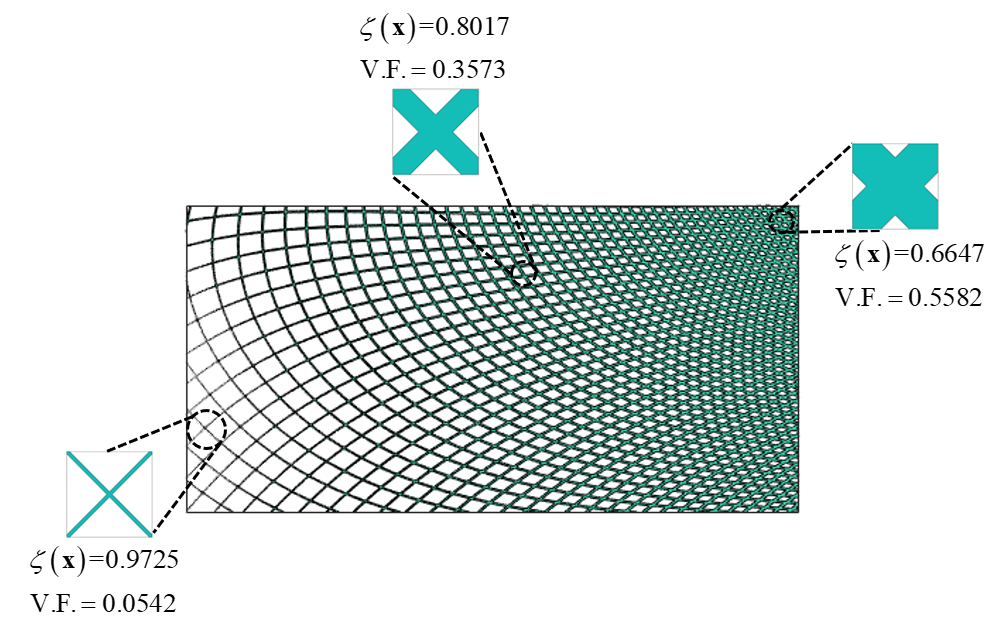}
    \caption{Optimised results for the uniform loading case shown in Fig.~\ref{f10}. Compared to the optimised configuration shown by Fig.~\ref{fig_optimised_uniform_load}, the volume fraction can vary in space now. Thus materials can be distributed better, so as to output a 26.7\% drop in the compliance value of the optimised GMC. Here ``V.F.'' stands for ``volume fraction''. \label{fig_optimised_vol_vary}}
\end{figure}

The study is also extended to cover the situations where the microscopic variables get activated, too. The loading condition here is the same as shown in Fig.~\ref{f10}. Then instead of a single matrix cell, a family of matrix cells, as exemplified in Fig.~\ref{fig_optimised_vol_vary}, are now in hand for GMC construction. The family of matrix cells is formed by allowing the thickness of the bars consisting the ``X''-shape configuration to change. Several issues making the present case more challenging than the previous one are worth being mentioned. First, an extra control parameter $\zeta(\mathbf{x})$ is introduced to identify the matrix cell of interest near the macroscopic position $\mathbf{x}$, and $\zeta$ also indicates the solid volume fraction there. Second, the TDF of use now is Eq.~\eqref{TDF_multiple_cells}. Third, by letting $d_1=\zeta$ in Eq.~\eqref{ML_target_function}, the microscopic design variables are also permitted to change. Now the number of the input arguments becomes four, meaning that the minimum amount of data for training should increase two. In this scenario, we following the procedure given in Sec.~\ref{Ch_3.2.1} to generate 40,000 sets of data, and a set of neural networks bearing an RMSE value no greater than 0.4\% are obtained through training. Finally for optimisation, the sensitivity of the microscopic design variable $\zeta$ is also needed, which can be calculated with a central difference scheme that is similar as Eq.~\eqref{central_diff}.

An optimised GMC can thus be worked out as shown in Fig.~\ref{fig_optimised_vol_vary}.
This time, the homogenised compliance of the optimised GMC is calculated to be 496, in comparison with 517 obtained with fine-mesh computation. Note that when the volume fraction is permitted to vary in space, the optimised compliance result records a 26.7\% drop if compared with the case with a single matrix cell as shown in Fig.~\ref{fig_optimised_uniform_load}(a). This is because materials are now allowed to be distributed in regions where the stress level is high, such as the top right corner in Fig.~\ref{fig_optimised_vol_vary}. Note that the volume fraction near the top right corner of the GMC shown in Fig.~\ref{fig_optimised_vol_vary} may reach 55\%, while the value is as low as 5\% near the bottom left corner.

\subsection{The case with a point load}

\begin{figure}[!ht]
    \centering
    \includegraphics[width=.5\textwidth]{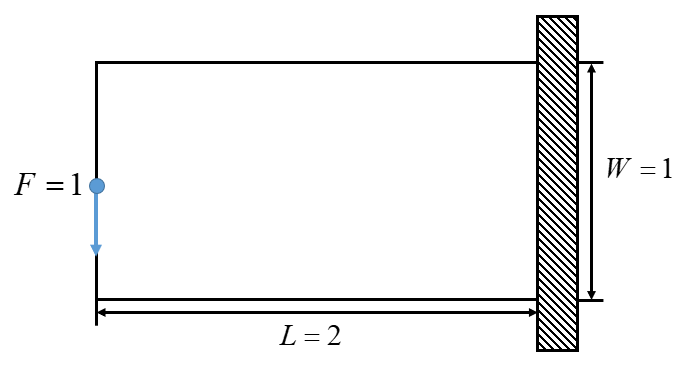}
    \caption{The case with a point load at the central point on the left edge\label{f13}}
\end{figure}

\begin{figure}[!ht]
    \centering
    \subfigure[]{\includegraphics[width=.5\textwidth]{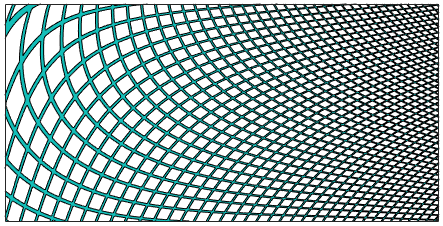}}
    \subfigure[]{\includegraphics[width=.35\textwidth]{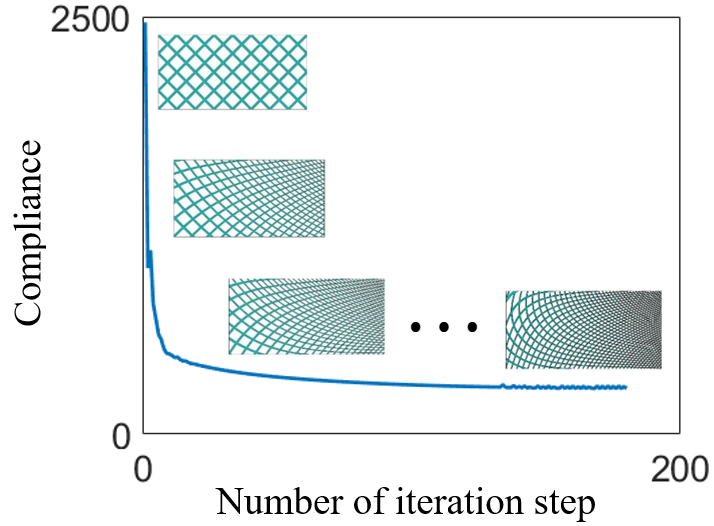}}
    \caption{Optimisation results for the case shown in Fig.~\ref{f13}. \label{fig_optimised_point_load}}
\end{figure}

Now we move to the case, as shown in Fig.~\ref{f13}, where the short beam is subject to a point load at the central point on the left edge.
Again we seek to find an optimised GMC which is filled by cells obtained by distorting a single matrix cell, as in the case of Fig.~\ref{fig_optimised_uniform_load}(a). The optimised GMC is shown in Fig.~\ref{fig_optimised_point_load}(a) with the optimisation procedure visualised in Fig.~\ref{fig_optimised_point_load}(b).
With the present homogenisation scheme, the compliance of the optimised GMC shown in Fig.~\ref{fig_optimised_point_load}(a) is calculated to be 252.3, which sees a 1.1\% deviation from the fine-mesh analysis over the same configuration. The result is also compared with the that obtained with the original zoning method \cite{xue2020speeding}, and a 3.3\% improvement in the optimised compliance value is seen.

\subsection{Computational efficiency\label{Sec_time}}
In the previous subsection, the accuracy shown by the present algorithm is demonstrated through comparisons with fine-mesh calculations. Here the computational efficiency delivered by the present scheme is shown. The issue will be addressed in two aspects: the clock time it takes to run optimisation, and how the performance is compared with that using the original zoning scheme \cite{xue2020speeding}.

On a standard desktop computer with one core activated, the clock time spent in running the optimisation algorithm is about 300 seconds for the two-dimensional examples presented here. Such a computational efficiency sits at a same level as the state-of-the art projection-based method \cite{groen2018homogenization}, qualifying the present scheme as one of the most efficient algorithms that can be used for stiffness optimisation of GMCs bearing complex microstructures.

More features about the computational efficiency from the method here can be found through a more detailed comparison with the performances of using the original zoning strategy \cite{xue2020speeding}. Here it should be noted that the original zoning strategy has its own meaning: it enables the realisation of the AHTO plus framework (at an acceptable cost). But the present scheme empowers the AHTO plus framework to a top level in a viewpoint of computational efficiency. The comparison results are summarised in Table~\ref{table_efficiency}, where the figures mean the clock time per iteration step during optimisation.
\begin{table}[!ht]
    \centering
    \renewcommand{\tablename}{Tabel}
    \caption{Computational efficiency in comparison with the original zoning strategy. Here ``OZS'' means the ``original zoning strategy''; ``N.A.'' stands for ``not available''. The data are measured in seconds.}
    \begin{tabular}{cccc}
        \toprule
        Number of zones & \makecell[c]{OZS \\non-parallel} & \makecell[c]{OZS \\in parallel}  & \makecell[c]{ Present \\ non-parallel}\\
        \midrule
        2 & 13.93 & 4.15 & 0.910 \\
        8 & 32.91 & 4.90 & 0.897 \\
        32 & 118.95 & 9.99 & 0.940 \\
        128 & 459.88 & 32.75 & 0.986 \\
        5000 & N.A. & N.A. & 3.436 \\		
        \bottomrule
    \end{tabular}
    \label{table_efficiency}	
\end{table}

The data are given where the design domains are partitioned into 2, 8, 32, 128 and 5000 zones. With a single core activated, it takes the original zoning strategy roughly 14 seconds to perform the optimisation for one step, which is roughly 14 times more than that from the present method. Moreover, the computational time increases with the number of subdomains for the original zoning strategy. This is because the microscopic cell problems have to be solved more times then. In contrast, the increase in the computational time with more subdomians is not significant in the present method.

Also it has been demonstrated that the efficiency of the zoning strategy may be enhanced with the activation of parallel computation. This can be observed from the third column of Table~\ref{table_efficiency}, where the data are obtained through parallel computing on a four-core processor. Nonetheless, the computational time for a single step is still several times more than that with the present method. Moreover, as the number of subdomains reaches 5,000, the original zoning strategy becomes disabled on the computational units mentioned above, because the memories are exhausted then. In contrast, the case with 5,000 subdomains is handled with ease on using the present scheme.

\section{Conclusion\label{Sec_conclusion}}

In this article, an optimisation scheme incorporating machine learning into the improved asymptotic-homogenisation-based topology optimisation framework \cite{zhu2019novel} is proposed, so as to accelerate the mechanical behaviour analysis and then the design of graded microstructural configurations. The introduction of the present algorithm addresses the issues on computational efficiency that limit the applicability of the AHTO plus method nowadays, where the microscopic cell problems have to be solved too many times. With neural networks trained against the data generated based on asymptotic analysis results, the characteristics of the unit cell problem can be taken into full account. Besides, the use of machine learning also enables the sensitivity calculations in a time-saving manner. A number of crucial issues on effectively integrating machine learning models with the AHTO plus formulation are discussed. For instance, the so-called low discrepancy sequence algorithm is employed so as to ensure that the generated data cover the entire design space more evenly. Besides, the Cholesky decomposition is also adopted so as to guarantee the positive definiteness of the homogenised elasticity tensor expressed by means of neural networks. More importantly, we also discuss the key issues on minimising the number of input arguments setting up a machine learning model, so as to bring down the requirements over the training data.

The accuracies and the efficiencies of the present scheme are demonstrated with numerical examples. For the issue of accuracy, the homogenised results are numerically shown to bear a deviation from the corresponding fine-mesh results no greater than 2\%. For the issue of efficiency, it takes roughly 300 seconds for the present algorithm to carry the optimisation out on a standard desktop computer, and this qualifies the AHTO plus method as one of the most efficient approaches for the compliance optimisation of configurations bearing complex microstructures. A comparative study over the computational efficiency against the original zoning strategy is also conducted, so as to reveal more features of the present scheme.

The AHTO plus method, now empowered by the use of machine learning, still needs to be improved in the following aspects.

Firstly, the representation of the mapping function $\mathbf{y}=\mathbf{y}(\mathbf{x})$ needs to be improved. This article continues the partition strategy, using polynomial functions as mapping functions. The disadvantage is that its design variables are the corresponding polynomial coefficients. Polynomials are in general unbounded, and the update of the polynomial coefficients has global relevance, that is, the local adjustment of the constituting cells is actually not permitted. This sometimes may lead to unacceptable distortion of matrix cells when the outputs of the mapping polynomials take quite large values. To solve this problem, one may consider using B-spline functions or other functions, where local adjustments of cell configurations are more flexible than with polynomial mapping functions.

Secondly, the machine learning model should be improved with a broader coverage. For instance, the matrix cells only take ``X''-shape configurations. The studies should be extended for various types of matrix cells, and sometimes their combinations. When the variety of matrix cell expands, the number of the input arguments setting up the machine learning models increases as well. This means one should further explore the interrelationships between these input arguments. For example, the chiral symmetry and rotational symmetry of the matrix cells should be formulated properly right before data collection.

Thirdly, there is a definite need for the present treatment to be generalised to cover three-dimensional examples. In three-dimensional scenarios, the GMC constituting cells are more likely to be ill-distorted \cite{xue2020generation}. This means that the optimisation then must be tuned with greater care. It can be expected that the success in developing a suitable algorithm for three-dimensional compliance optimisation of GMC is built on a solid understanding over the two-dimensional situations, especially over the issues raised above.

\section*{Acknowledgement}
The financial supports from the National Natural Science Foundation of China (11772076, 11732004, 11821202, 12172074) are gratefully acknowledged.

\section*{Appendix 1}
Here we show that the homogenised elasticity tensor $\mathbb{C}^{ \mathrm{H}}$ stay in variant, as all the entries of the Jacobian matrix $\mathbf{J}$ get multiplied or divided by a same coefficient.

Suppose $J_{i j}=\alpha J_{i j}^{\prime}$, where $\alpha$ can be any constant. Substituting it into governing equation of the cell problem, i.e., Eq.~\eqref{eq8}, we obtain
\begin{equation}\label{eq20}
    \alpha J_{m j}^{\prime} \frac{\partial}{\partial \bar{Y}_{m}}\left(\overline{\mathbb{C}}_{i j k l} \alpha J_{n l}^{\prime} \frac{\partial \xi_{k}^{s t}}{\partial \bar{Y}_{n}}\right)=\alpha J_{m j}^{\prime} \frac{\partial \tilde{\mathrm{C}}_{i j s t}}{\partial \bar{Y}_{m}}, \quad \text { in } \Upsilon_{\mathrm{F}},
\end{equation}
which becomes
\begin{equation}\label{eq21}
    J_{m j}^{\prime} \frac{\partial}{\partial \bar{Y}_{m}}\left(\overline{\mathbb{C}}_{i j k} J_{n l}^{\prime} \frac{\partial\left(\alpha \xi_{k}^{s t}\right)}{\partial \bar{Y}_{n}}\right)=J_{m j}^{\prime} \frac{\partial \tilde{\mathbb{C}}_{i s t}}{\partial \bar{Y}_{m}}, \quad \text { in } \Upsilon_{\text {P }}.
\end{equation}
Eq.~\eqref{eq21} is effectively the governing equation for the third-order $\xi_{k}^{s t^{\prime}}=\alpha \xi_{k}^{s t}$ tensor corresponding to $\mathbf{J}^{\prime}$. And substituting the homogenisation modulus calculation formula~\eqref{eq7}, we derive
\begin{equation}\label{eq22}
    \begin{aligned}
        \mathbb{C}_{i j k l}^{\mathrm{H}^{\prime}} &=\mathbb{C}_{i j k k} \cdot\left|\Upsilon_{\mathrm{P}}^{s}\right|-\mathbb{C}_{i j s t} J_{n t}^{\prime} \int_{\Upsilon_{\mathrm{P}}} \frac{\partial \xi_{s}^{k \prime}}{\partial \bar{Y}_{n}} \mathrm{~d} \overline{\mathbf{Y}} \\
        &=\mathbb{C}_{i j k l} \cdot\left|\Upsilon_{\mathrm{P}}^{s}\right|-\mathbb{C}_{i j s t} \frac{J_{n t}}{\alpha} \int_{\Upsilon_{\mathrm{P}}} \frac{\partial\left(\alpha \xi_{s}^{k t}\right)}{\partial \bar{Y}_{n}} \mathrm{~d} \overline{\mathbf{Y}} \\
        &=\mathbb{C}_{i j k l} \cdot\left|\Upsilon_{\mathrm{P}}^{s}\right|-\mathbb{C}_{i j s t} J_{n t} \int_{\Upsilon_{\mathrm{P}}} \frac{\partial \xi_{s}^{k d}}{\partial \bar{Y}_{n}} \mathrm{~d} \overline{\mathbf{Y}} \\
        &=\mathbb{C}_{i j k l}^{\mathrm{H}} \quad, \quad i, j, k, l=1, \ldots, N.
    \end{aligned}
\end{equation}

\section*{References}
\bibliography{reference}
\end{document}